\documentstyle[12pt]{article}
\begin{document}
\newtheorem{thm}{Theorem}
\newtheorem{prop}{Proposition}
\newtheorem{cor}{Corollary}
\newtheorem{ex}{Example}
\newtheorem{exs}{Examples}
\newtheorem{lem}{Lemma}
\newtheorem{rem}{Remark}
\newcommand{\bt}{\begin{thm}}
\newcommand{\et}{\end{thm}}
\newcommand{\bl}{\begin{lem}}
\newcommand{\el}{\end{lem}}
\newcommand{\bp}{\begin{prop}}
\newcommand{\ep}{\end{prop}}
\newcommand{\bc}{\begin{cor}}
\newcommand{\ec}{\end{cor}}
\newcommand{\p}{{{\bf Proof.\,\,}}}
\def\NN{I\!\! N}
\def\RR{I\!\! R}
\def\HH{I\!\! H}
\def\QQ{I\!\!\!\! Q}
\def\CC{I\!\!\!\! C}
\def\ZZ{Z\!\!\! Z}

\newcommand{\mod}[3]{\mbox{${#1} \equiv{#2}\bmod{#3}$}}
\newcommand{\G} {{\bf G}}
\newcommand{\biindice}[3]%
{
\renewcommand{\arraystretch}{0.5}
\begin{array}[t]{c}
#1\\
{\scriptstyle #2}\\
{\scriptstyle #3}
\end{array}
\renewcommand{\arraystretch}{1}
}
\newcommand{\chit}{\mbox{${\tilde\chi}$}}
\def\qed{\hfill \vrule height4pt width4pt depth2pt}
\baselineskip6mm \centerline{{\Large Commutativity conditions on derivations and Lie ideals}}
\centerline{{\Large   in  $\sigma$-prime rings}} \vspace*{0,6cm} \centerline{{\large L.
Oukhtite, S. Salhi \& L. Taoufiq}} \centerline{\small{Universit\'e
Moulay Isma\"\i l, Facult\'e des Sciences et Techniques}}
\centerline{\small{D\'epartement de Math\'ematiques, Groupe
d'Alg\`ebre et Applications}} \centerline{\small{ B. P. 509
Boutalamine, Errachidia; Maroc}}
\centerline{\small{oukhtitel@hotmail.com, salhi@fastmail.fm,
lahcentaoufiq@yahoo.com}}
\begin{abstract}
\vspace*{-0,2cm} \noindent Let $R$ be a 2-torsion free
$\sigma$-prime ring, $U$ a nonzero square closed $\sigma$-Lie ideal
of $R$ and let $d$ be a derivation of $R.$ In this paper it is shown
that:\\
1) If $d$ is centralizing on $U,$ then $d=0$ or $U\subseteq Z(R).$\\
2) If either $d([x,y])=0$ for all $x,y\in U,$ or $[d(x),d(y)]=0$ for
all $x,y \in U$ and $d$ commutes with $\sigma$ on $U,$ then $d=0$ or
$U\subseteq Z(R).$
\end{abstract}

\section{ Introduction}
\vspace*{-0,5cm} Throughout this paper, $R$ will represent an
associative ring with center $Z(R).$ Recall that $R$ is said to be
$2$-torsion free if whenever $2x=0,$ with $x\in R,$ then $x=0.$ $R$
is prime if $aRb=0$ implies that $a=0$ or $b=0$ for all $a$ and $b$
in $R.$ If $\sigma$ is an involution in $R,$ then $R$ is said  to be
$\sigma$-prime if $a R b =a R\sigma (b)=0$ implies that  $a=0$ or
$b=0$. It is obvious that every prime ring equipped with an
involution $\sigma$ is also $\sigma$-prime, but the converse need
not be true in general.  An additive mapping $d : R\rightarrow R$ is
said to be a derivation if $d(xy)=d(x)y+xd(y)$ for all $x,y$ in $R.$
A mapping $F : R\rightarrow R$ is said to be centralizing on a
subset $S$ of $R$ if $[F(s),s]\in Z(R)$ for all $s\in S.$ In
particular, if $[F(s),s]=0$ for all $s\in S$, then $F$ is commuting
on $S.$ In all that follows $Sa_{\sigma}(R)$ will denote the set of
symmetric and skew-symmetric elements of $R$; i.e.,
$Sa_{\sigma}(R)=\{x\in R/ \sigma(x)=\pm x\}.$ For any $x,y \in R,$
the commutator $xy-yx$ will be denoted by $[x,y].$ An additive subgroup
$U$ of $R$ is said to be a Lie ideal of $R$ if $[u,r]\in U$ for all
$u\in U$ and $r\in R.$ A Lie ideal $U$ which satisfies
$\sigma(U)\subseteq U$ is called a $\sigma$-Lie ideal. If $U$ is a
Lie (resp. $\sigma$-Lie) ideal of $R,$ then $U$ is called a square
closed Lie (resp. $\sigma$-Lie) ideal if $u^{2}\in U$ for all $u\in
U.$ Since $(u+v)^{2}\in U$ and $[u,v]\in U$, we see that $2uv\in U$
for all $u,v \in U.$ Therefore, for all $r\in R$ we get
$2r[u,v]=2[u,rv]-2[u,r]v \in U$ and $2[u,v]r=2[u,vr]-2v[u,r]\in U,$
so that $2R[U,U]\subseteq U$ and $2[U,U]R\subseteq U.$ This remark
will be freely used in the whole paper.\\Many works concerning the
relationship between commutativity of a ring and the behavior of
derivations defined on this ring have been studied. The first
important result in this subject is Posner's Theorem,  which states
that the existence of a nonzero centralizing derivation on a prime
ring forces this ring to be commutative
($\cite{Pos}$). This result has been generalized by many authors in several ways.\\
In $\cite{Her}$, I. N. Herstein proved that if $R$ is a prime ring
of characteristic not $2$ which has a nonzero derivation $d$ such
that $[d(x), d(y)] = 0$ for all $x, y \in R,$ then  $R$ is
commutative. Motivated by this result, H. E. Bell, in $\cite{Bell}$,
studied derivations $d$ satisfying $d([x,y])=0$ for all $ x, y \in
R$. In $\cite{O.S2}$ and $\cite{O.S5}$, L. Oukhtite
and S. Salhi generalized these results to $\sigma$-prime rings. In
particular, they proved that if $R$ is a $2$-torsion free
$\sigma$-prime ring equipped with a  nonzero derivation which is
 centralizing on $R,$ then $R$ is necessarily
commutative.\\Our purpose in this paper is to extend these results
to square closed $\sigma$-Lie ideals. \vspace*{-0.5cm}
 \section{Preliminaries and results}
 In order to prove our main theorems, we shall need the following lemmas.
\vspace*{-0,5cm}
\bl \label{l1}

\em{(\cite{adv}, Lemma 4)}  If $U\not\subset Z(R)$ is a $\sigma$-Lie
ideal of a $2$-torsion free $\sigma$-prime ring $R$ and $a,b\in R$
such that $aUb=\sigma(a)Ub=0$ or $aUb=aU\sigma(b)=0,$ then $a=0$ or $b=0.$

  \el
\vspace*{-0.6cm}
\bl\label{l2}

\em{(\cite{IJA.1.1 2007}, Lemma 2.3)}  Let $0\neq U$ be a
$\sigma$-Lie ideal of a $2$-torsion free $\sigma$-prime ring $R.$ If
$[U,U]= 0,$ then $U \subseteq Z(R).$

\el
\vspace*{-0.7cm}
\bl \label{l3}

\em{(\cite{O.S4}, Lemma 2.2)} Let $R$ be a $2$-torsion free
$\sigma$-prime ring and $U$ a nonzero $\sigma$-Lie ideal of $R.$ If
$d$ is a derivation of $R$ which commutes with $\sigma$ and
satisfies $d(U)=0,$ then either $d=0$ or $U\subseteq Z(R).$

\el
\vspace*{-0.4cm}
{\bf Remark.}
One can easily verify that  Lemma 3 is still valid if  the condition that $d$ commutes with $\sigma$ is replaced by $d\circ \sigma= - \sigma \circ d.$
\vspace*{-0.4cm}
\bt\label{t1}

Let $R$ be a 2-torsion free $\sigma$-prime ring and $U$  a square
closed $\sigma$-Lie ideal of $R.$ If d is a derivation of $R$
satisfying $[d(u), u] \in Z(R)$ for all $u\in U,$  then $U \subseteq
Z(R)$ or $d=0$.

\et
 \p Suppose that $\;U\not\subseteq Z(R).$ As $\;[d(x),x] \in Z(R)\;$ for all $x\in U,$
by linearization $[d(x),y]+[d(y),x]\in Z(R) $ for all $x,y\in U.$
Since char$R\neq 2,$ the fact that $[d(x),x^2]+[d(x^2),x] \in Z(R) $
yields
 $x[d(x),x] \in Z(R)$ for all $x\in U;$  hence $$[r,x][d(x),x]=0\;\;\mbox{
for all }\;\; x\in U, r\in R, $$ and therefore $[d(x),x]^{2}=0$ for
all $ x\in U.$ Since $[d(x),x] \in Z(R)$,
$$[d(x),x]R[d(x),x]\sigma([d(x),x])=0\;\;\mbox{for all}\;\;x\in U$$  and the $\sigma$-primeness of
$R$ yields $\;[d(x),x]=0\;$ or $\;[d(x),x]\sigma([d(x),x])=0.$ If
$[d(x),x]\sigma([d(x),x])=0$, then $[d(x),x]R\sigma([d(x),x])=0;$ and
the fact that $[d(x),x]^{2}=0$ gives
$$ [d(x),x]R\sigma([d(x),x])=[d(x),x]R[d(x),x]=0.$$Since $R$ is $\sigma$-prime, we obtain
 $$[d(x),x]=0\;\mbox{ for all }\;\; x\in U.$$ Let us
consider the map $\delta : R\longmapsto R$ defined by
$\delta(x)=d(x) + \sigma\circ d\circ \sigma(x).$ One can easily
verify that $\delta$ is a derivation of $R$ which commutes with
$\sigma$ and satisfies
$$ [\delta(x),x]=0\;\mbox{ for all }\;\; x\in U.$$ Linearizing this
equality, we obtain
$$[\delta(x),y]+[\delta(y),x]=0\;\;\mbox{ for all}\;\;x,y\in U.$$ Writing
$2xz$ instead of $y$ and using char$R\neq 2,$ we find that
$$\delta(x)[x,z]=0\;\;\mbox{ for all}\;\;x,z\in U.$$ Replacing $z$ by $2zy$ in this
equality,  we conclude that $\delta(x)z[x,y]=0,$ so that
\begin{equation} \label{equ1}
\delta(x)U[x,y]=0\;\;\mbox{for all}\;\;x,y\in U.
\end{equation} By virtue of Lemma 1, it then follows
that
$$
 \delta(x)=0\;\;\mbox{or}\;\;[x,U]=0,\;\;\mbox{for all}\;\;x\in U\cap Sa_{\sigma}(R).
$$
 Let $u\in U$. Since $u-\sigma(u)\in U\cap Sa_{\sigma}(R),$ it follows that $$\delta(u-\sigma(u))=0\;\;\mbox{or}\;\;
[u-\sigma(u),U]=0.$$ If $\delta(u-\sigma(u))=0,$ then $\delta(u)\in
Sa_{\sigma}(R)$ and $(1)$ yields $\;\delta(u)=0\;$ or $\;[u,U]=0.$
If $[u-\sigma(u),U]=0,$ then $ [u,y]=[\sigma(u),y]$ for all $y\in U$
and $(1)$ assures that  $$\delta(u)U [u,y]= 0=\delta(u)U
 \sigma([u,y]),\;\;\mbox{for all}\;\;y\in U.$$ Applying Lemma 1, we find that $\delta(u)=0$ or
$[u,U]=0.$  Hence, $U$ is a union of two additive subgroups $G_1$
and $G_2$, where $$G_{1}=\{ u\in U\;\mbox{such
that}\;\delta(u)=0\}\;\;\mbox{and}\;\;G_{2}=\{ u\in U\;\mbox{such
that}\;[u,U]=0\}.$$ Since a group cannot be a union of two of its
proper subgroups, we are forced to  $U =G_{1}$ or $U=G_{2}.$ Since
$U\not\subseteq Z(R),$ Lemma 2 assures that $U =G_{1}$ and therefore
$\delta(U)=0.$ Now applying Lemma $3$, we get $\delta=0$ and
therefore  $d\circ \sigma = - \sigma\circ d.$ As $[d(x),x]=0$ for
all $\; x\in U,$ in view of the above Remark, similar reasoning
leads to $d=0.$ \qed
 \vspace*{-0.3cm}
  \bc\label{c1}

\em{(\cite{O.S5}, Theorem 1.1)} Let $R$ be a 2-torsion free
$\sigma$-prime ring and $d$ a nonzero derivation of $R.$ If d is
centralizing on $R,$ then $R$ is commutative.

\ec
 \vspace*{-0.6cm}
\bt\label{t2}

Let $U$ be a  square closed $\sigma$-Lie ideal of a  $2$-torsion
free $\sigma$-prime ring $R$ and $d$ a derivation of $R$ which
commutes with $\sigma$ on $U.$\\If  $[d(x),d(y)]=d([y,x])$ for all
$x,y \in U$, then  $d=0$ or  $ U\subseteq Z(R).$

\et
\vspace*{-0,5cm}
\p Suppose that  $ U\not\subset Z(R).$ We have
\begin{equation} \label{equ1}[d(x),d(y)]=d([y,x])\;\;\mbox{for all}\;\;x,y\in U.
\end{equation}
Substituting $2xy$ for $y$ in $(2)$ and using char$R\neq 2,$ we get
\begin{equation} \label{equ1}d(x)[y,x]=[d(x),x]d(y)+d(x)[d(x),y]\;\;\mbox{for all}\;\;x,y\in U.
\end{equation}
Replacing $y$ by $2[y,z]x$ and using (3), we find that
\begin{equation} \label{equ1}[d(x),x][y,z]d(x)+d(x)[y,z][d(x),x]=0\;\;\mbox{for all}\;\;x,y, z\in U.
\end{equation}
Replace $y$ by $2[y,z]d(x)$ in $(3)$  to get
\begin{equation} \label{equ1}d(x)[y,z][d(x),x]-[d(x),x][y,z]d^2(x)=0\;\;\mbox{for all}\;\;x,y, z\in U.
\end{equation}
From $(4)$ and $(5)$ we obtain
\begin{equation} \label{equ1}[d(x),x][y,z](d(x)+d^2(x))=0\;\;\mbox{for all}\;\;x,y,z\in U.
\end{equation}
Writing $2[u,v](d(x)+d^2(x))y$ instead of $y$  in $(6),$ where $ u,v\in U,$ we obtain
$[d(x),x][u,v]z(d(x)+d^2(x))y(d(x)+d^2(x))=0,$ so  that
\begin{equation} \label{equ1}[d(x),x][u,v]z(d(x)+d^2(x))U(d(x)+d^2(x))=0\;\;\mbox{for all}\;\;x, u,v, z\in U.
\end{equation}
If $x\in U\cap Sa_{\sigma}(R)$, then  Lemma 1 together with $(7)$
assures that
$$d(x)+d^2(x)=0\;\;\mbox{ or }\;\;
[d(x),x][u,v]z(d(x)+d^2(x))=0\;\;\mbox{ for all}\;\;u,v,z\in
U.$$Suppose that  $[d(x),x][u,v]z(d(x)+d^2(x))=0.$ Then
\begin{equation} [d(x),x][u,v]U(d(x)+d^2(x))=0.
\end{equation}
Since $d$ commutes with $\sigma$ and $x\in Sa_{\sigma}(R),$ in view
of $(8)$, Lemma 1 gives
\begin{equation} d(x)+d^2(x)=0\;\;\mbox{or}\;\;[d(x),x][u,v]=0\;\;\mbox{for all}\;\;
u,v\in U.
\end{equation}
If $[d(x),x][u,v]=0$, then replacing $u$ by $2uw$ in $(9)$ where
$w\in U,$ we obtain
\begin{equation} [d(x),x]U[u,v]=0.
\end{equation}
As $\;\sigma(U)=U$ and $\;[U,U]\neq 0,$ by $(10)$, Lemma 2 yields
that $\;[d(x),x]=0$. Thus, in any event,
$$
\mbox{either}\;\;[d(x),x]=0\;\;\mbox{or}\;\;d(x)+d^2(x)=0 \;
\;\mbox{for all}\;\;x\in U\cap Sa_{\sigma}(R).
$$
Let $x\in U.$   Since $x+\sigma(x)\in U\cap Sa_{\sigma}(R),$ either $d(x+\sigma(x)) + d^2( x+\sigma(x))=
0$ or $[d(x+\sigma(x)), x+\sigma(x)]=0.$\\
If $d(x+\sigma(x)) + d^2( x+\sigma(x))=0,$ then $d(x)+ d^2( x) \in
Sa_{\sigma}(R)$   and $(7)$ yields that $d(x)+ d^2( x)=0$ or
$[d(x),x][u,v]U(d(x)+d^2(x))=0.$\\If
$[d(x),x][u,v]U(d(x)+d^2(x))=0,$ once again using $d(x)+ d^2( x) \in
Sa_{\sigma}(R),$ we find that $d(x)+ d^2( x) =0,$ or
$[d(x),x][u,v]$ for all $u,v\in U,$ in which case $[d(x),x]=0.$\\Now
suppose that  $[d(x+\sigma(x)), x+\sigma(x)]=0.$ As $x-\sigma(x)\in
U\cap Sa_{\sigma}(R)
$  we have to distinguish two cases:\\
1) If $d(x-\sigma(x))+ d^2(x-\sigma(x))=0,$ then $d(x)+d^2(x)\in Sa_{\sigma}(R)$. Reasoning as above we get $d(x)+
 d^2( x) =0$ or $[d(x),x]=0.$\\
2) If  $[d(x-\sigma(x)), x-\sigma(x)]=0$, then $[d(x),x]\in Sa_{\sigma}(R).$
Replace $u$ by $2yu$ in $(7)$, with $y\in U$ , to get
$[d(x),x]y[u,v]z(d(x)+d^2(x))U(d(x)+d^2(x)) =0,$ so that
\begin{equation} \label{equ1}[d(x),x]U[u,v]z(d(x)+d^2(x))U(d(x)+d^2(x)) =0\;\;\mbox{for all}\;\;x,u,v, z\in U.
\end{equation}
Since $[d(x),x]\in Sa_{\sigma}(R),$ from $(11)$ it follows that
$$[d(x),x]=0\;\;\mbox{or}\;\;[u,v]U(d(x)+d^2(x))U(d(x)+d^2(x))
=0\;\mbox{ for all }\;u,v\in U.$$Suppose
$[u,v]U(d(x)+d^2(x))U(d(x)+d^2(x)) =0.$ As $\sigma(U)=U$ and
$[U,U]\neq 0,$ then
\begin{equation} \label{equ1}(d(x)+d^2(x))U(d(x)+d^2(x)) =0.
\end{equation}
In $(6)$, write $2[u,v](d(x)+d^2(x))r$ instead of $z$, where $u,v\in
U$ and $r\in R,$  to obtain
 \begin{equation} \label{equ1}
 [d(x),x][u,v]y(d(x)+d^2(x))r(d(x)+d^2(x)) =0,\;\mbox{for
 all}\;u,v,y\in U, r\in R. \end{equation}
Replacing $r$ by $r\sigma(d(x)+d^2(x))z$ in $(13)$, where $z\in U,$
we find that
$$[d(x),x][u,v]y(d(x)+d^2(x))r\sigma(d(x)+d^2(x))z(d(x)+d^2(x)) =0,$$ which leads us to
\begin{equation} \label{equ1}[d(x),x][u,v]y(d(x)+d^2(x))U\sigma(d(x)+d^2(x))U(d(x)+d^2(x)) =0.
\end{equation}Since $\sigma(d(x)+d^2(x))U(d(x)+d^2(x))$ is invariant under $\sigma$, by virtue of $(14),
$ Lemma 1 yields
$$\sigma(d(x)+d^2(x))U(d(x)+d^2(x))
=0\;\;\mbox{or}\;\;[d(x),x][u,v]y(d(x)+d^2(x))=0.$$If
$\sigma(d(x)+d^2(x))U(d(x)+d^2(x)) =0,$ then $(12)$ implies that
$d(x)+d^2(x)=0.$ \\Now assume that
\begin{equation} \label{equ1}
[d(x),x][u,v]y(d(x)+d^2(x))=0\;\;\mbox{for all}\;\;u,v,y\in U.
\end{equation}
Replace $v$ by $2wv$ in $(15)$, where $w\in U$, and use $(15)$ to
get$$[d(x),x]w[u,v]y(d(x)+d^2(x))=0,$$ so that
\begin{equation}
[d(x),x]U[u,v]y(d(x)+d^2(x))=0\;\mbox{for all}\;\;u,v,y\in U.
\end{equation} As $[d(x),x]\in Sa_{\sigma}(R),$ $(16)$ yields $[u,v]U(d(x)+d^2(x))=0,$  in which case
$d(x)+d^2(x)=0,$ or $[d(x),x]=0$.\\In conclusion, for all $x\in U$ we have either $[d(x),x]=0$ or $d(x)+d^2(x)=0.$\\
Now let $x\in U$ such that $d(x)+d^2(x)=0.$  In $(2)$, put $y=2[y,z]d(x)$ to get
\begin{equation} \label{equ1}d([y,z])[d(x),x] -[[y,z],x] d(x) +[d(x),[y,z]]d(x)=[y,z][d(x),x].
\end{equation}
If in $(2)$ we put $y=2[y,z]x,$ we get
\begin{equation} \label{equ1}[[y,z],x]d(x)  = [d(x),[y,z]] d(x) + d([y,z])[d(x),x]=0.
\end{equation}
From $(17)$ and $(18)$ it then follows that
$$[y,z][d(x),x]=0\;\;\mbox{for all}\;\;y,z\in U ,$$ hence $\;[y,z]U [d(x),x]
=0$ for all $y,z\in U.$ Applying Lemma 1, this  leads to
$$[d(x),x]=0,\;\;\mbox{for all}\;\;x\in U.$$ By virtue of Theorem 1,
this yields that $d=0.$ \qed\\
\\
Note that if $d$ is a derivation of $R$ which acts as an
anti-homomorphism on $U$, then $d$ satisfies the condition
$[d(x),d(y)]=d([y,x])$ for all $x,y \in U.$ Thus we have the
following corollary.
  \vspace*{-0,4cm}
\bc

\em{(\cite{O.S4}, Theorem 1.1)} Let $d$ be a derivation of a
$2$-torsion free $\sigma$-prime ring $R$ which acts as an
anti-homomorphism on a nonzero square closed $\sigma$-Lie ideal $U$
of $R.$ If $d$ commutes with $\sigma$, then either $d=0$ or $
U\subseteq Z(R).$

\ec
  \vspace*{-0,6cm}
  \bt\label{t3}

Let $U$ be a square closed $\sigma$-Lie ideal of a  $2$-torsion free
$\sigma$-prime ring $R$ and $d$  a derivation of $R.$  If either
$d([x,y])=0$ for all $x,y\in U,$ or  $[d(x),d(y)]=0$ for all $x,y
\in U$ and $d$ commutes with $\sigma$ on $U,$  then  $d=0$ or $
U\subseteq Z(R).$

 \et
\vspace*{-0,5cm}
 \p
Suppose that $ U\not\subseteq Z(R).$ Assume that $d([x, y])=0\;$ for
all $x,y \in U .$ Let $\delta$ be  the derivation of $ R$ defined by
$\delta(x) = d(x) +\sigma\circ d \circ \sigma(x).$\\Clearly,
$\delta$ commutes with $\sigma$ and $\delta([x,y])=0$ for all
$x,y\in U,$ so that
\begin{equation}
\label{equ1}[\delta(x),y]= [\delta(y),x]\;\;\mbox{ for all} \;\;x,y
\in U.
\end{equation} Writing $[x,y]$ instead of $y$ in
$(19),$ we find that
\begin{equation} \label{equ1}
[\delta(x), [x,y]]=0 \;\;\mbox{ for all} \;\; x,y \in U.
\end{equation} Replacing $x$ by $x^2$ in $(19)$, we
conclude that
\begin{equation} \label{equ1}
\delta(x)[x,y]+[x,y]\delta(x)=0 \;\;\mbox{ for all} \;\; x,y \in U.
\end{equation}
As char$R\neq 2$, from $(20)$ and $(21)$ it follows that
\begin{equation} \label{equ1} \delta(x)[x,y]=0
\;\;\mbox{ for all} \;\; x,y \in U.
\end{equation}
Replacing $y$ by $2zy$ in $(22),$ we get  $\delta(x)z[x,y]=0 ,$ so
that
$$ \delta(x)U[x,y]=0  \;\;\mbox{ for all} \;\;
x,y \in U.$$ From the proof of Theorem 1, we conclude that
$\delta=0$ and thus $d\circ \sigma= -\sigma \circ d.$ Since $d$
satisfies $d( [x,y])=0$ for all $x,y\in U,$ by similar reasoning,
we are forced to $d=0.$\\Now assume that $d$ commutes with $\sigma$
and satisfies $[d(x),d(y)]=0$ for all $x,y \in U.$ The fact that
$[d(x),d(2xy)]=0 $  implies that
\begin{equation} \label{equ1}
d(x)[d(x),y]+[d(x),x]d(y)=0\;\;\mbox{ for all} \;\;x,y \in U .
\end{equation}
Replace $y$ by $2[y,z]d(u)$ in $(23)$, where $z,u\in U,$  to find
that
\begin{equation} \label{equ1}
[d(x),x][y,z]d^{2}(u)=0\;\;\mbox{ for all}
\;\;x,y, u \in U .
\end{equation}
Write $2[s,t]d^2(w)y$ instead of $y$ in $(24)$, where $s,t,w\in U,$
thereby concluding that $[d(x),x]z[s,t]d^2(w)yd^2(u)=0.$ Accordingly,
\begin{equation} \label{equ1}
[d(x),x]z[s,t]d^2(w)Ud^2(u)=0\;\;\mbox{for all}\;\;s,t, u,w,x\in
U.
\end{equation} Since $d$ commutes with $\sigma$
and $\sigma(U)=U,$ using $(25)$ we find that
$$d^2(U)=0\;\;\mbox{or}\;\; [d(x),x]U[s,t]d^2(w)=0. $$ Suppose that
\begin{equation} \label{equ1}
[d(x),x]U[s,t]d^2(w)=0\;\;\mbox{for all}\;\;s,t, w,x\in
U.
\end{equation} Replacing $t$ by $2tv$ in $(26)$, where $v\in
U,$ we are forced to $$[d(x),x][s,t]vd^2(w)=0$$ and hence
\begin{equation} \label{equ1} [d(x),x][s,t]Ud^2(w)=0\;\;\mbox{for
all}\;\;s,t, w,x\in U.
\end{equation} Since $\sigma(U)=U$ and $d$ commutes with $\sigma$,
then $(27)$ implies that either $d^2(U)=0,$ or $[d(x),x][s,t]=0$ for all
$s,t,x\in U,$ in which case $[d(x),x]=0$ for all $x\in U.$\\Thus, in any
event, we find that
$$d^2(U)=0\;\;\mbox{or}\;\;[d(x),x]=0\;\;\mbox{for all}\;x\in U.$$
If $d^2(U)=0,$ then [\cite{IJA.1.1 2007}, Theorem 1.1]  assures that
$d=0.$\\If $[d(x),x]=0$ for all $x\in U$, then Theorem 1 yields
$d=0.$ \qed
 \vspace*{-0,3cm}
\bc

\em{(\cite{O.S2}, Theorem 3.3)} Let $d$ be a nonzero derivation of a
$2$-torsion free $\sigma$-prime ring $R.$ If $d([x,y])=0$ for all
$x,y\in R,$ then $R$ is commutative.

\ec

\end{document}